\documentclass[11pt,a4paper]{article}
\pdfoutput=1

\usepackage[utf8]{inputenc}
\usepackage{amsmath, amssymb, amsthm}
\usepackage{geometry}
\usepackage{hyperref}
\usepackage{graphicx}
\usepackage{tikz}
\usepackage{booktabs}
\usepackage{tabularx}
\usepackage{caption}
\usepackage{float}
\usetikzlibrary{shapes.geometric, arrows.meta, positioning, calc}

\title{\textbf{Nineteen to the Dozen: Embedding the Neo-Riemannian Tonnetz into a Cyclic $19_3$ Symmetric Configuration}}
\author{
    \textbf{Paweł Nurowski} \\
    \small Centrum Fizyki Teoretycznej, Polska Akademia Nauk, Al. Lotników 32/46, 02-668 Warszawa, Poland \\
    \small Guangdong Technion Israel Institute of Technology, No. 241, Daxue Road, Shantou, China \\
    \small \texttt{nurowski@cft.edu.pl}
}
\date{\today}

\newtheorem{theorem}{Theorem}
\newtheorem{definition}{Definition}
\newtheorem{observation}{Observation}

\begin{document}

\maketitle
\begin{abstract}
This paper bridges combinatorial geometry, music theory, and biomechanics to solve the fundamental challenge of embedding classical Western harmony into the microtonal 19-tone equal temperament (19-TET). Inspired by Roger Penrose's observations on the mathematical elegance of 19-TET, we provide the theoretical foundation for a physical 19-TET acoustic piano currently under construction. However, playing classical 12-TET music on such an instrument poses a topological problem: emvedding the classical Euler-Riemann Tonnetz into the 19-TET universe inevitably distorts structural chords, creating dissonant ``wolves.'' By formalizing these harmonic spaces as incidence configurations (the $12_3$ and $19_3$ graphs) and utilizing integer cuts in our optimization model, we exhaustively prove that exactly 32 out of 36 Neo-Riemannian harmonic connections can be preserved. We demonstrate a strict 5-fold degeneracy of this optimum: there exist exactly 5 mathematically equivalent local packings for the wolf chords. Among these, we identify a unique canonical realization in which the 14 excised vertices form a perfectly contiguous geometric void along the primary Hamiltonian cycle. We reveal that the 4 inevitably broken edges represent the exact topological scars of the historical enharmonic diesis, and we formulate the Vicentino Hypothesis regarding 16th-century microtonal composition. Finally, to make this theoretical geometry physically playable, we design a novel 19-TET split-key keyboard, formalized through a biomechanical cost function that optimizes the performer's hand span. This work provides the complete theoretical, historical, and ergonomic blueprint for the next generation of microtonal acoustic instruments.
\end{abstract}

\section{Introduction: From Decaphony to Nineteen}

The genesis of this research stems from a direct conversation with Sir Roger Penrose following the premiere of a decaphonic (10-TET) acoustic piano. The instrument was designed and built as a collaborative effort between the author, Aleksander Bogucki, Andrzej Włodarczyk, and the team of renowned jazz pianist Leszek Możdżer. Upon hearing of the successful construction of the decaphonic piano, Penrose remarked: \textit{``What a shame! It should have been 19!''} When asked why, he provided his 2000 paper, \textit{The Heritage of Pythagoras: Nineteen to the Dozen} \cite{penrose2000}. This serendipitous exchange served as the author's first exposure to the profound mathematical and acoustic elegance of the 19-tone equal temperament (19-TET) and directly inspired the search for topological connections between the classical 12-TET and 19-TET spaces.

The interplay between combinatorics and musical harmony has recently entered a rigorous geometric phase. In their recent work, Boland and Hughston \cite{boland2025} demonstrated that the classical Euler-Riemann Tonnetz can be formalized not merely as a heuristic network, but as a strict incidence geometry. Specifically, by representing the 12 major chords as points and the 12 minor chords as lines in a Euclidean plane, they established that the Tonnetz is isomorphic to the self-dual $\{12_3\}$ combinatorial configuration known as the $D_{222}$ of Daublebsky von Sterneck \cite{daublebsky1895}. 

However, this formalization raises an equally profound algebraic question: is the $D_{222}$ topology an exclusive property of Western acoustics, or is it merely one solution within a broader mathematical landscape? In our previous investigations of 12-TET and 10-TET spaces \cite{nurowski12, nurowski10}, we treated harmonic systems purely as solutions to the linear Diophantine system $t+s \equiv q \pmod n$, revealing that the specific topological structure of standard Western harmony is one of exactly twelve mathematically isomorphic `shadow' systems in 12-TET.

Building upon these foundations and Penrose's suggestion, the present paper solves a specific combinatorial problem: the geometric embedding of the 12-TET Neo-Riemannian Tonnetz into the microtonal, prime-numbered universe of 19-TET. Because 19 is a prime number, the underlying topology cannot be decomposed into a sub-grid. To achieve this, we first establish the rigorous rules for generating the Tonnetze in both spaces, define the subset of 19-TET that approximates Just Intonation, and define the optimization constraints required to mend the topological graph ``scars'' caused by the enharmonic diesis. Ultimately, we point out the specific, severe construction challenges that engineers will face when building a true, acoustic piano based on these theories.

\section{Formal Definitions of the Harmonic Networks}

\subsection{The $12_3$ Tonnetz ($D_{222}$)}
The Neo-Riemannian Tonnetz in 12-TET consists of 12 Major triads and 12 minor triads \cite{cohn2012}. Mathematically, it forms the $12_3$ symmetric configuration cataloged historically as $D_{222}$ \cite{daublebsky1895, alazemi2014}.

\begin{definition}[12-TET Tonnetz]
Let $\mathbb{Z}_{12}$ represent the pitch classes of 12-TET. The vertices of the Tonnetz graph are defined as:
\begin{itemize}
    \item Major triads: $M^{12}_k = \{k, (k+4) \bmod 12, (k+7) \bmod 12\}$
    \item Minor triads: $m^{12}_k = \{k, (k+3) \bmod 12, (k+7) \bmod 12\}$
\end{itemize}
for $k \in \{0, \dots, 11\}$. Rather than defining edges via abstract set intersection, we define them explicitly via voice-leading transformations. A Major triad $M^{12}_k$ connects to exactly three minor triads: $m^{12}_k$ (Parallel), $m^{12}_{k+4}$ (Relative), and $m^{12}_{k+9}$ (Leading-Tone). Reciprocally, a minor triad $m^{12}_k$ connects to three Major triads: $M^{12}_k$, $M^{12}_{k+8}$, and $M^{12}_{k+3}$.
\end{definition}

\subsection{The $19_3$ Microtonal Configuration}

In 19-TET, the intervals mapping to the major and minor thirds are respectively 6 and 5 steps, while the perfect fifth is 11 steps.

\begin{definition}[19-TET Harmonic Network]
Let $\mathbb{Z}_{19}$ represent the pitch classes of 19-TET. The vertices are defined as:
\begin{itemize}
    \item Major triads: $M^{19}_r = \{r, (r+6) \bmod 19, (r+11) \bmod 19\}$
    \item Minor triads: $m^{19}_r = \{r, (r+5) \bmod 19, (r+11) \bmod 19\}$
\end{itemize}
for $r \in \{0, \dots, 18\}$. Just as in 12-TET, rather than defining edges via abstract set intersection, we define them explicitly via voice-leading transformations. A Major triad $M^{19}_r$ connects to exactly three minor triads: $m^{19}_r$ (Parallel), $m^{19}_{r+6}$ (Relative), and $m^{19}_{r+14}$ (Leading-Tone). Reciprocally, a minor triad $m^{19}_r$ connects to three Major triads: $M^{19}_r$, $M^{19}_{r+13}$, and $M^{19}_{r+5}$.
\end{definition}

While the total number of generic configurations grows hyper-exponentially \cite{betten2000} achieving the number of over 7.6 billion $19_3$ configurations \cite{grunbaum2009}, we restrict our attention strictly to \textit{cyclic} configurations.

\begin{observation}
The fact that there exist exactly three non-isomorphic cyclic $19_3$ configurations is a direct consequence of difference set enumeration. A cyclic $19_3$ configuration is generated by a base block $B = \{0, a, b\} \pmod{19}$. The configuration's incidence structure is entirely determined by its difference set $\Delta(B) = \{\pm a, \pm b, \pm(b-a)\} \pmod{19}$. For the Levi graph to be a valid configuration (girth $\ge 6$), the three positive differences $x, y, z \in \{1, \dots, 9\}$ extracted from $\Delta(B)$ must be strictly distinct and satisfy either $x+y=z$ or $x+y+z=19$.

Exhaustive enumeration yields exactly 16 unique difference sets satisfying $x+y=z$ and exactly 5 sets satisfying $x+y+z=19$, resulting in 21 valid difference sets. Two cyclic configurations generated by blocks $B_1$ and $B_2$ are isomorphic if and only if their difference sets belong to the same orbit under the action of the affine multiplier group. 

Let us define this group action explicitly. The multiplicative group $\mathbb{Z}_{19}^*$ consists of 18 elements $\{1, 2, \dots, 18\}$. Because the edges of our geometric configuration are undirected, a step of size $x$ is geometrically equivalent to a step of size $-x \pmod{19}$ (for instance, a distance of $1$ is equivalent to a distance of $18$). We formalize this symmetry by taking the quotient group $G = \mathbb{Z}_{19}^* / \{1, -1\}$. This quotient operation identifies each multiplier $m$ with its additive inverse $19-m$, effectively pairing the 18 elements into exactly 9 equivalence classes: $\{1, 18\}, \{2, 17\}, \dots, \{9, 10\}$. Therefore, $G$ is a group of order 9. Furthermore, this group is cyclic and is generated by the multiplier $m=2$, since successive powers of 2 modulo 19 (when equating $x$ with $-x$) systematically cycle through all 9 classes.

The action of this 9-element group $G$ (via scalar multiplication modulo 19) on the 21 valid difference sets partitions them into precisely 3 orbits: two regular orbits of length 9, and one short orbit of length 3 (stabilized by the subgroup generated by $m=7$). Thus, there are precisely three affine equivalence classes of cyclic $19_3$ configurations.
\end{observation}

From these three, we uniquely select the configuration generated by the cyclic offset $\{0, 6, 14\}$. This choice is strictly dictated by the acoustic reality of 19-TET: because the best approximation of a Major Third is exactly 6 steps, and a minor third is 5 steps, preserving the common-tone voice-leading rules of the Tonnetz dictates that a Major triad $M^{19}_r$ must connect to exactly three minor triads: $m^{19}_r$ (sharing the 5th), $m^{19}_{r+6}$ (sharing the minor 3rd), and $m^{19}_{r+14}$ (sharing the Major 3rd, since $14 \equiv -5 \pmod{19}$).

We place its 38 vertices on a Hamiltonian cycle alternating steps of $+14$ and $-6$ (the 19-TET equivalent of the Circle of Fifths).

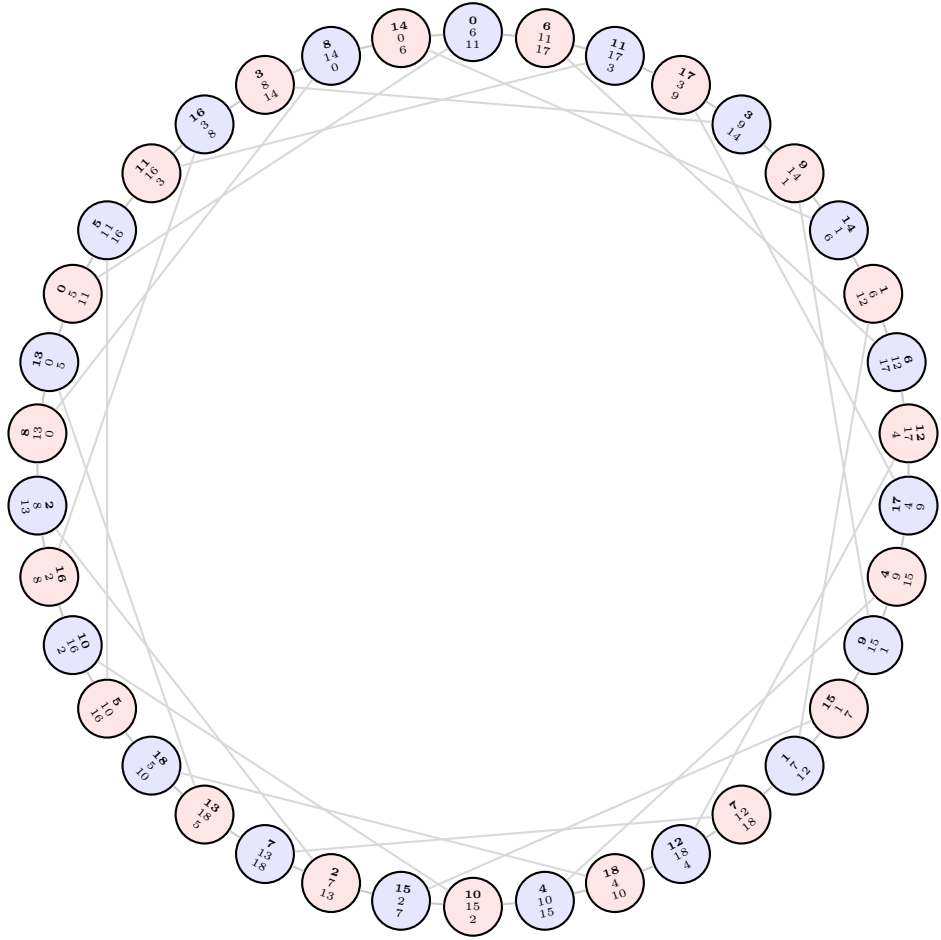
\begin{figure}[H]
\centering
\begin{tikzpicture}[scale=0.85, transform shape]
    \def\N{19}
    \def\R{6.8cm}
    
    \foreach \val [count=\i from 0] in {0, 8, 16, 5, 13, 2, 10, 18, 7, 15, 4, 12, 1, 9, 17, 6, 14, 3, 11} {
        \pgfmathsetmacro{\angle}{360/(2*\N) * (2*\i) + 90}
        \pgfmathsetmacro{\vtwo}{int(mod(\val+6, 19))}
        \pgfmathsetmacro{\vthree}{int(mod(\val+11, 19))}
        \pgfmathsetmacro{\angmod}{mod(\angle + 3600, 360)}
        \pgfmathsetmacro{\rot}{\angmod > 180 ? \angmod + 90 : \angmod - 90}
        
        \node[circle, draw, thick, fill=blue!10, minimum size=0.9cm, inner sep=0pt] (M\val) at (\angle: \R) {};
        \node[font=\fontsize{5}{6}\selectfont, rotate=\rot, align=center] at (\angle: \R) {\textbf{\val} \\[-0.3ex] \vtwo \\[-0.3ex] \vthree};
    }
    
    \foreach \val [count=\i from 0] in {14, 3, 11, 0, 8, 16, 5, 13, 2, 10, 18, 7, 15, 4, 12, 1, 9, 17, 6} {
        \pgfmathsetmacro{\angle}{360/(2*\N) * (2*\i + 1) + 90}
        \pgfmathsetmacro{\vtwo}{int(mod(\val+5, 19))}
        \pgfmathsetmacro{\vthree}{int(mod(\val+11, 19))}
        \pgfmathsetmacro{\angmod}{mod(\angle + 3600, 360)}
        \pgfmathsetmacro{\rot}{\angmod > 180 ? \angmod + 90 : \angmod - 90}
        
        \node[circle, draw, thick, fill=red!10, minimum size=0.9cm, inner sep=0pt] (m\val) at (\angle: \R) {};
        \node[font=\fontsize{5}{6}\selectfont, rotate=\rot, align=center] at (\angle: \R) {\textbf{\val} \\[-0.3ex] \vtwo \\[-0.3ex] \vthree};
    }
    
    \foreach \i in {0,...,18} {
        \pgfmathparse{int(mod(\i+6,19))} \let\nSix\pgfmathresult
        \draw[thick, gray!40] (M\i) -- (m\nSix);
        \pgfmathparse{int(mod(\i+14,19))} \let\nFourteen\pgfmathresult
        \draw[thick, gray!40] (M\i) -- (m\nFourteen);
        \draw[thick, gray!30] (M\i) -- (m\i);
    }
\end{tikzpicture}
\caption{The pristine $19_3$ configuration showing full connectivity and radial pitch-class triples according to the $\{0, 6, 14\}$ offset.}
\label{fig:19_3_full}
\end{figure}

\section{The Embedding Rules: Just Intonation and Wolf Chords}

To embed the 24 vertices of $12_3$ into the 38 vertices of $19_3$, we must define a mapping rule. We select a specific subset of 12 tones within 19-TET:
$$S_{JI} = \{0, 2, 3, 5, 6, 8, 10, 11, 13, 14, 16, 17\} \subset \mathbb{Z}_{19}$$
This subset is chosen because these 12 pitches provide the absolute best mathematical approximations to 5-limit Just Intonation within the 19-TET scale. We define the natural embedding map $\phi: \mathbb{Z}_{12} \hookrightarrow \mathbb{Z}_{19}$ into this subset:
$$ \phi = \{0{\mapsto}0, 1{\mapsto}2, 2{\mapsto}3, 3{\mapsto}5, 4{\mapsto}6, 5{\mapsto}8, 6{\mapsto}10, 7{\mapsto}11, 8{\mapsto}13, 9{\mapsto}14, 10{\mapsto}16, 11{\mapsto}17\} $$

\subsection{The Emergence of Wolf Chords}
When the 24 classical triads of $12_3$ are passed through $\phi$, a structural dichotomy emerges:
\begin{enumerate}
    \item \textbf{16 Perfect Chords:} Exactly 16 chords (8 Major, 8 minor) map to valid $19_3$ vertices. Their internal intervals perfectly match the $\{6, 5\}$ requirement of 19-TET triads. For example, C-Major $M_0^{12} = \{0, 4, 7\} \xrightarrow{\phi} \{0, 6, 11\} = M_0^{19}$. These 16 vertices become rigid geometric anchors in our graph.
    \item \textbf{8 Wolf Chords:} The remaining 8 chords suffer severe interval distortion because their mapped intervals cross the "missing" microtonal steps of our 12-tone subset. For example, D-Major $M_2^{12} = \{2, 6, 9\} \xrightarrow{\phi} \{3, 10, 14\}$. The intervals here are $10-3=7$ and $14-10=4$, which blatantly violates the $\{6,5\}$ structure of a legal 19-TET triad. These are the ``wolf chords.''
\end{enumerate}

\section{MILP Optimization and Topological Legalization}

A naive embedding is impossible because the 8 wolf chords do not correspond to any valid nodes in the $19_3$ space. To resolve this, we formulated a Maximum Common Edge Subgraph (MCES) problem using a Mixed Integer Linear Programming (MILP) solver.

\noindent
\textbf{Optimization Rules:} 
\begin{enumerate}
\item The 16 perfect chords are strictly anchored to their exact mathematical coordinates in the $19_3$ Levi graph.
\item 2. The solver must place the 8 wolf chords by assigning each of them a unique position among the 22 remaining available vertices in $19_3$. In this way, we populate exactly the $16+8=24$ vertices of the target $D_{222}$ graph.
\item \textbf{Objective Function:} Maximize the number of preserved edges $e_{ij}$ (corresponding to Neo-Riemannian P, L, R voice-leading transformations) from the original $D_{222}$ topology.
\end{enumerate}

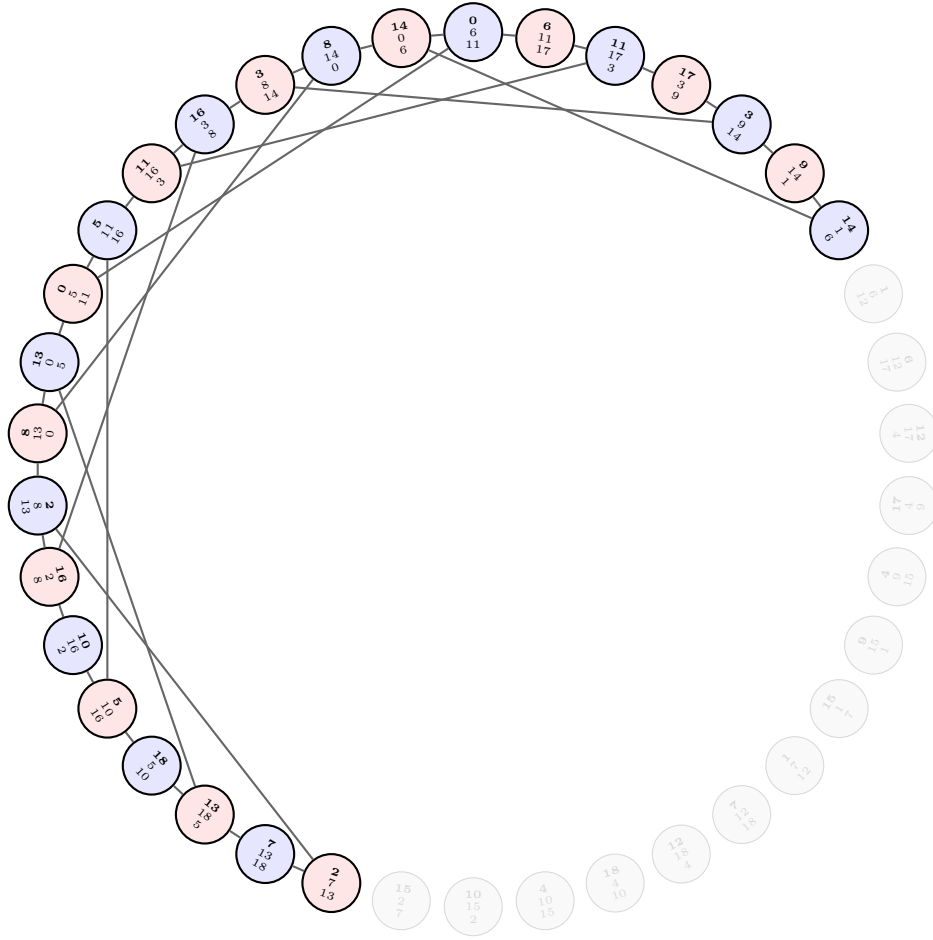
\begin{figure}[H]
\centering
\begin{tikzpicture}[scale=0.85, transform shape]
    \def\N{19}
    \def\R{6.8cm}
    
    \foreach \val [count=\i from 0] in {0, 8, 16, 5, 13, 2, 10, 18, 7, 15, 4, 12, 1, 9, 17, 6, 14, 3, 11} {
        \pgfmathsetmacro{\angle}{360/(2*\N) * (2*\i) + 90}
        \pgfmathsetmacro{\vtwo}{int(mod(\val+6, 19))}
        \pgfmathsetmacro{\vthree}{int(mod(\val+11, 19))}
        \pgfmathsetmacro{\angmod}{mod(\angle + 3600, 360)}
      
        \pgfmathsetmacro{\rot}{\angmod > 180 ? \angmod + 90 : \angmod - 90}
        
        \ifnum\i>8 \ifnum\i<16
            \node[circle, draw=gray!30, fill=gray!5, minimum size=0.9cm, inner sep=0pt] (M\val) at (\angle: \R) {};
            \node[gray!30, font=\fontsize{5}{6}\selectfont, rotate=\rot, align=center] at (\angle: \R) {\textbf{\val} \\[-0.3ex] \vtwo \\[-0.3ex] \vthree};
        \else
            \node[circle, draw, thick, fill=blue!10, minimum size=0.9cm, inner sep=0pt] (M\val) at (\angle: \R) {};
            \node[font=\fontsize{5}{6}\selectfont, rotate=\rot, align=center] at (\angle: \R) {\textbf{\val} \\[-0.3ex] \vtwo \\[-0.3ex] \vthree};
        \fi \else
            \node[circle, draw, thick, fill=blue!10, minimum size=0.9cm, inner sep=0pt] (M\val) at (\angle: \R) {};
            \node[font=\fontsize{5}{6}\selectfont, rotate=\rot, align=center] at (\angle: \R) {\textbf{\val} \\[-0.3ex] \vtwo \\[-0.3ex] \vthree};
        \fi
    }
    
    \foreach \val [count=\i from 0] in {14, 3, 11, 0, 8, 16, 5, 13, 2, 10, 18, 7, 15, 4, 12, 1, 9, 17, 6} {
        \pgfmathsetmacro{\angle}{360/(2*\N) * (2*\i + 1) + 90}
        \pgfmathsetmacro{\vtwo}{int(mod(\val+5, 19))}
        \pgfmathsetmacro{\vthree}{int(mod(\val+11, 19))}
        \pgfmathsetmacro{\angmod}{mod(\angle + 3600, 360)}
        \pgfmathsetmacro{\rot}{\angmod > 
        180 ? \angmod + 90 : \angmod - 90}
        
        \ifnum\i>8 \ifnum\i<16
            \node[circle, draw=gray!30, fill=gray!5, minimum size=0.9cm, inner sep=0pt] (m\val) at (\angle: \R) {};
            \node[gray!30, font=\fontsize{5}{6}\selectfont, rotate=\rot, align=center] at (\angle: \R) {\textbf{\val} \\[-0.3ex] \vtwo \\[-0.3ex] \vthree};
        \else
            \node[circle, draw, thick, fill=red!10, minimum size=0.9cm, inner sep=0pt] (m\val) at (\angle: \R) {};
            \node[font=\fontsize{5}{6}\selectfont, rotate=\rot, align=center] at (\angle: \R) {\textbf{\val} \\[-0.3ex] \vtwo \\[-0.3ex] \vthree};
        \fi \else
            \node[circle, draw, thick, fill=red!10, minimum size=0.9cm, inner sep=0pt] (m\val) at (\angle: \R) {};
            \node[font=\fontsize{5}{6}\selectfont, rotate=\rot, align=center] at (\angle: \R) {\textbf{\val} \\[-0.3ex] \vtwo \\[-0.3ex] \vthree};
        \fi
    }
    
    \draw[thick, black!60] (M0) -- (m0); \draw[thick, black!60] (M0) -- (m6); \draw[thick, black!60] (M0) -- (m14);
    \draw[thick, black!60] (M2) -- (m2); \draw[thick, black!60] (M2) -- (m8); \draw[thick, black!60] (M2) -- (m16);
    \draw[thick, black!60] (M3) -- (m3); \draw[thick, black!60] (M3) -- (m9); \draw[thick, black!60] (M3) -- (m17);
    \draw[thick, black!60] (M5) -- (m5); \draw[thick, black!60] (M5) -- (m11); \draw[thick, black!60] (M5) -- (m0);
    \draw[thick, black!60] (M7) -- (m13); \draw[thick, black!60] (M7) -- (m2);
    \draw[thick, black!60] (M8) -- (m8); \draw[thick, black!60] (M8) -- (m14); \draw[thick, black!60] (M8) -- (m3);
    \draw[thick, black!60] (M10) -- (m16); \draw[thick, black!60] (M10) -- (m5);
    \draw[thick, black!60] (M11) -- (m11); \draw[thick, black!60] (M11) -- (m17); \draw[thick, black!60] (M11) -- (m6);
    \draw[thick, black!60] (M13) -- (m13); \draw[thick, black!60] (M13) -- (m0); \draw[thick, black!60] (M13) -- (m8);
    \draw[thick, black!60] (M14) -- (m14); \draw[thick, black!60] (M14) -- (m9);
    \draw[thick, black!60] (M16) -- (m16); \draw[thick, black!60] (M16) -- (m3); \draw[thick, black!60] (M16) -- (m11);
    \draw[thick, black!60] (M18) -- (m5); \draw[thick, black!60] (M18) -- (m13);
\end{tikzpicture}
\caption{The $19_3$ configuration during the excision process (representing the canonical Solution 1). The 14 unmapped elements are grayed out, leaving an unclosed contiguous void.}
\label{fig:19_3_excision}
\end{figure}

When executing this embedding, an extraordinary geometric phenomenon occurs. For the primary optimal configuration found by the solver (which we designate as the canonical Solution 1), the 14 unmapped (deleted) vertices form a single, perfectly contiguous block along the Hamiltonian cycle. Figure \ref{fig:19_3_excision} captures this system \textit{during the excision phase}, delineating this continuous 14-node chunk as desaturated ``ghost'' vertices.

The MILP solver successfully mends this gap by matching the 8 non-matching ``Wolf'' chords from the naive mapping onto legal structural coordinates, sacrificing exactly 4 edges to secure a maximum overall edge count of 32 out of 36. Table \ref{tab:milp_subs_canonical} documents these exact substitutions for this canonical realization.
\begin{table}[H]
\centering
\caption{MILP Topological Repair: Resolving the 8 Wolf Chords (Canonical Solution 1)}
\label{tab:milp_subs_canonical}
{\tiny \begin{tabularx}{\textwidth}{@{}l X X X@{}}
\toprule
\textbf{12-TET Orig.} & \textbf{Naive Mapping via $\phi$} & \textbf{MILP Repaired Node} & \textbf{Tonnetz Role} \\
\midrule
$M^{12}_{2}$ (D)   & $\{3, 10, 14\}$ (Invalid) & $M^{19}_{3} = \{3, 9, 14\} \in 19_3$ & Restores 2 edges, breaks 1 (Scar) \\&&&\\
$M^{12}_{4}$ (E)   & $\{6, 13, 17\}$ (Invalid) & $M^{19}_{7} = \{7, 13, 18\} \in 19_3$ & Restores 2 edges, breaks 1 (Scar) \\&&&\\
$M^{12}_{9}$ (A)   & $\{14, 2, 6\}$ (Invalid)  & $M^{19}_{14} = \{14, 1, 6\} \in 19_3$ & Restores 2 edges, breaks 1 (Scar) \\&&&\\
$M^{12}_{11}$ (B)  & $\{17, 5, 10\}$ (Invalid) & $M^{19}_{18} = \{18, 5, 10\} \in 19_3$ & Restores 2 edges, breaks 1 (Scar) \\&&&\\
$m^{12}_{1}$ (C$\sharp$m)  & $\{2, 6, 13\}$ (Invalid) & $m^{19}_{2} = \{2, 7, 13\} \in 19_3$ & Bridges local gap \\&&&\\
$m^{12}_{6}$ (F$\sharp$m)  & $\{10, 14, 2\}$ (Invalid) & $m^{19}_{9} = \{9, 14, 1\} \in 19_3$ & Bridges local gap \\&&&\\
$m^{12}_{8}$ (G$\sharp$m)  & $\{13, 17, 5\}$ (Invalid) & $m^{19}_{13} = \{13, 18, 5\} \in 19_3$ & Bridges local gap \\&&&\\
$m^{12}_{11}$ (Bm) & $\{17, 3, 10\}$ (Invalid) & $m^{19}_{17} = \{17, 3, 9\} \in 19_3$ & Bridges local gap \\
\bottomrule
\end{tabularx}}
\end{table}

\begin{figure}[H]
\centering
\definecolor{mygold}{RGB}{218,165,32}
\begin{minipage}{0.48\textwidth}
\centering
\begin{tikzpicture}[scale=0.8, transform shape]
    \def\N{19}
    \def\R{5.3cm}
    
    \foreach \val/\ismilp/\isactive [count=\i from 0] in {
        0/0/1, 8/0/1, 16/0/1, 5/0/1, 13/0/1, 2/0/1, 10/0/1, 18/1/1, 7/1/1, 
        15/0/0, 4/0/0, 12/0/0, 1/0/0, 9/0/0, 17/0/0, 6/0/0, 
        14/1/1, 3/1/1, 11/0/1%
    } {
        \ifnum\isactive=1
            \pgfmathsetmacro{\angle}{360/(2*\N) * (2*\i) + 90}
            \pgfmathsetmacro{\vtwo}{int(mod(\val+6, 19))}
            \pgfmathsetmacro{\vthree}{int(mod(\val+11, 19))}
            \pgfmathsetmacro{\angmod}{mod(\angle + 3600, 360)}
            \pgfmathsetmacro{\rot}{\angmod > 180 ? \angmod + 90 : \angmod - 90}
            
            \ifnum\ismilp=1
                \node[circle, draw=mygold, line width=1.6pt, fill=blue!10, minimum size=0.9cm, inner sep=0pt] (M\val) at (\angle: \R) {};
            \else
                \node[circle, draw=black, thick, fill=blue!10, minimum size=0.9cm, inner sep=0pt] (M\val) at (\angle: \R) {};
            \fi
            \node[font=\fontsize{5}{6}\selectfont, rotate=\rot, align=center] at (\angle: \R) {\textbf{\val} \\[-0.3ex] \vtwo \\[-0.3ex] \vthree};
        \fi
    }
    
    \foreach \val/\ismilp/\isactive [count=\i from 0] in {
        14/0/1, 3/0/1, 11/0/1, 0/0/1, 8/0/1, 16/0/1, 5/0/1, 13/1/1, 2/1/1, 
        10/0/0, 18/0/0, 7/0/0, 15/0/0, 4/0/0, 12/0/0, 1/0/0, 
        9/1/1, 17/1/1, 6/0/1%
    } {
        \ifnum\isactive=1
            \pgfmathsetmacro{\angle}{360/(2*\N) * (2*\i + 1) + 90}
            \pgfmathsetmacro{\vtwo}{int(mod(\val+5, 19))}
            \pgfmathsetmacro{\vthree}{int(mod(\val+11, 19))}
            \pgfmathsetmacro{\angmod}{mod(\angle + 3600, 360)}
            \pgfmathsetmacro{\rot}{\angmod > 180 ? \angmod + 90 : \angmod - 90}
            
            \ifnum\ismilp=1
                \node[circle, draw=mygold, line width=1.6pt, fill=red!10, minimum size=0.9cm, inner sep=0pt] (m\val) at (\angle: \R) {};
            \else
                \node[circle, draw=black, thick, fill=red!10, minimum size=0.9cm, inner sep=0pt] (m\val) at (\angle: \R) {};
            \fi
            \node[font=\fontsize{5}{6}\selectfont, rotate=\rot, align=center] at (\angle: \R) {\textbf{\val} \\[-0.3ex] \vtwo \\[-0.3ex] \vthree};
        \fi
    }
    
    \draw[thick, black] (M0) -- (m0); \draw[thick, black] (M0) -- (m6); \draw[thick, black] (M0) -- (m14);
    \draw[thick, black] (M2) -- (m2); \draw[thick, black] (M2) -- (m8); \draw[thick, black] (M2) -- (m16);
    \draw[thick, black] (M3) -- (m3); \draw[thick, black] (M3) -- (m9); \draw[thick, black] (M3) -- (m17);
    \draw[thick, black] (M5) -- (m5); \draw[thick, black] (M5) -- (m11); \draw[thick, black] (M5) -- (m0);
    \draw[thick, black] (M7) -- (m13); \draw[thick, black] (M7) -- (m2);
    \draw[thick, black] (M8) -- (m8); \draw[thick, black] (M8) -- (m14); \draw[thick, black] (M8) -- (m3);
    \draw[thick, black] (M10) -- (m16); \draw[thick, black] (M10) -- (m5);
    \draw[thick, black] (M11) -- (m11); \draw[thick, black] (M11) -- (m17); \draw[thick, black] (M11) -- (m6);
    \draw[thick, black] (M13) -- (m13); \draw[thick, black] (M13) -- (m0); \draw[thick, black] (M13) -- (m8);
    \draw[thick, black] (M14) -- (m14); \draw[thick, black] (M14) -- (m9);
    \draw[thick, black] (M16) -- (m16); \draw[thick, black] (M16) -- (m3); \draw[thick, black] (M16) -- (m11);
    \draw[thick, black] (M18) -- (m5); \draw[thick, black] (M18) -- (m13);
    
    \draw[thick, dashed, mygold, line width=2pt] (M7) -- (m6);
    \draw[thick, dashed, mygold, line width=2pt] (M10) -- (m9);
    \draw[thick, dashed, mygold, line width=2pt] (M14) -- (m2);
    \draw[thick, dashed, mygold, line width=2pt] (M18) -- (m17);
\end{tikzpicture}
\end{minipage}\hfill
\begin{minipage}{0.48\textwidth}
\centering
\begin{tikzpicture}[scale=0.7, transform shape]
    \def\R{5.3cm}
    
    \foreach \name/\val/\vtwo/\vthree/\ismilp/\ismajor [count=\i from 0] in {
        M0/0/6/11/0/1, m14/14/0/6/0/0, M8/8/14/0/0/1, m3/3/8/14/0/0,
        M16/16/3/8/0/1, m11/11/16/3/0/0, M5/5/11/16/0/1, m0/0/5/11/0/0,
        M13/13/0/5/0/1, m8/8/13/0/0/0, M2/2/8/13/0/1, m16/16/2/8/0/0,
        M10/10/16/2/0/1, m5/5/10/16/0/0, M18/18/5/10/1/1, m13/13/18/5/1/0,
        M7/7/13/18/1/1, m2/2/7/13/1/0, M14/14/1/6/1/1, m9/9/14/1/1/0,
        M3/3/9/14/1/1, m17/17/3/9/1/0, M11/11/17/3/0/1, m6/6/11/17/0/0%
    } {
        \pgfmathsetmacro{\angle}{90 - 360/24 * \i}
        \pgfmathsetmacro{\angmod}{mod(\angle + 3600, 360)}
        \pgfmathsetmacro{\rot}{\angmod > 180 ? \angmod + 90 : \angmod - 90}
        
        \ifnum\ismajor=1
            \ifnum\ismilp=1
                \node[circle, draw=mygold, line width=1.6pt, fill=blue!10, minimum size=0.9cm, inner sep=0pt] (\name) at (\angle: \R) {};
            \else
                \node[circle, draw=black, thick, fill=blue!10, minimum size=0.9cm, inner sep=0pt] (\name) at (\angle: \R) {};
            \fi
        \else
            \ifnum\ismilp=1
                \node[circle, draw=mygold, line width=1.6pt, fill=red!10, minimum size=0.9cm, inner sep=0pt] (\name) at (\angle: \R) {};
            \else
                \node[circle, draw=black, thick, fill=red!10, minimum size=0.9cm, inner sep=0pt] (\name) at (\angle: \R) {};
            \fi
        \fi
        \node[font=\fontsize{5}{6}\selectfont, rotate=\rot, align=center] at (\angle: \R) {\textbf{\val} \\[-0.3ex] \vtwo \\[-0.3ex] \vthree};
    }
    
    \draw[thick, black!80] (M0) -- (m14); \draw[thick, black!80] (m14) -- (M8);
    \draw[thick, black!80] (M8) -- (m3); \draw[thick, black!80] (m3) -- (M16);
    \draw[thick, black!80] (M16) -- (m11); \draw[thick, black!80] (m11) -- (M5);
    \draw[thick, black!80] (M5) -- (m0); \draw[thick, black!80] (m0) -- (M13);
    \draw[thick, black!80] (M13) -- (m8); \draw[thick, black!80] (m8) -- (M2);
    \draw[thick, black!80] (M2) -- (m16); \draw[thick, black!80] (m16) -- (M10);
    \draw[thick, black!80] (M10) -- (m5); \draw[thick, black!80] (m5) -- (M18);
    \draw[thick, black!80] (M18) -- (m13); \draw[thick, black!80] (m13) -- (M7);
    \draw[thick, black!80] (M7) -- (m2); \draw[thick, dashed, mygold, line width=2pt] (m2) -- (M14);
    \draw[thick, black!80] (M14) -- (m9); \draw[thick, black!80] (m9) -- (M3);
    \draw[thick, black!80] (M3) -- (m17); \draw[thick, black!80] (m17) -- (M11);
    \draw[thick, black!80] (M11) -- (m6); \draw[thick, black!80] (m6) -- (M0);
    
    \draw[thick, black!30] (M0) -- (m0); \draw[thick, black!30] (m14) -- (M14);
    \draw[thick, black!30] (M8) -- (m8); \draw[thick, black!30] (m3) -- (M3);
    \draw[thick, black!30] (M16) -- (m16); \draw[thick, black!30] (m11) -- (M11);
    \draw[thick, black!30] (M5) -- (m5); \draw[thick, black!30] (M13) -- (m13);
    \draw[thick, black!30] (M2) -- (m2);
    
    \draw[thick, dashed, mygold, line width=2pt] (M10) -- (m9);
    \draw[thick, dashed, mygold, line width=2pt] (M18) -- (m17);
    \draw[thick, dashed, mygold, line width=2pt] (M7) -- (m6);
\end{tikzpicture}
\end{minipage}
\caption{Dual representation of the legalized $D_{222}$ sub-graph for canonical Solution 1. Left: Placed within the original 19-TET static embedding space. Right: Unraveled completely along its 24-node structural Hamiltonian cycle. All 8 MILP substituted nodes and the 4 repaired harmonic seams are explicitly highlighted in gold, verifying the exact $i \leftrightarrow i+7$ Moebius-Kantor topology.}
\label{fig:side_by_side}
\end{figure}
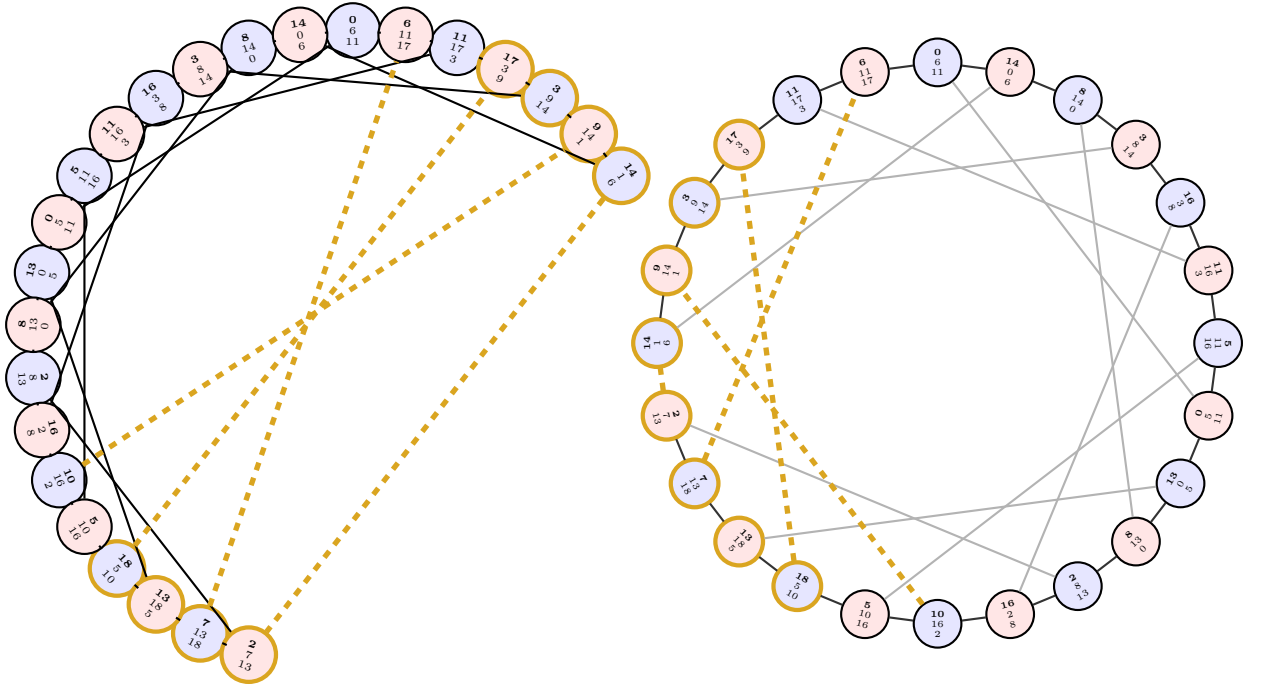

Following this table, Figure \ref{fig:side_by_side} provides a comprehensive dual visualization for the canonical Solution 1. On the left, the repaired $D_{222}$ space is displayed maintaining the original 19-TET embedding coordinates. On the right, the network is subjected to a structural Hamiltonian unraveling, revealing the true topological connectivity as a continuous ring with symmetric chordal crossings mapping strictly to $i \to i+7$.

\subsection{Musical and Acoustic Perspective: Approach for Performers and Tuners}
For a musician, piano tuner, or musicologist, the operation of the optimization algorithm is not an abstract numerical procedure, but a mathematical reflection of tangible, auditory work on the instrument. When naive mapping $\phi$ projects a classical D-Major chord ($M^{12}_2$) into the 19-TET space, its pure, consonant intervals (the major and minor thirds, which possess ideal widths of 6 and 5 microtonal steps in 19-TET, respectively) are severely distorted, resulting in intervals that are 7 and 4 steps wide. Such a chord becomes a harsh, out-of-tune ``wolf.''

At this stage, the algorithm performs a task identical to that of a piano tuner correcting an instrument by ear, or the intuition of an avant-garde Renaissance composer like Nicola Vicentino \cite{vicentino1555}: it microtonally ``retunes'' individual pitches of the wolf chord, searching among the remaining 22 available, perfectly tuned chords of the 19-TET system for the one that lies ``closest.'' However, the criterion for this closeness is strictly functional—the solver selects a replacement chord that maximizes common tones with its natural neighbors within the harmonic web (the smoothness of voice-leading). Consequently, although the absolute pitches are slightly shifted, the pianist's fingers on the keyboard can still navigate an intact, fluid grid of 32 out of 36 original pathways, allowing the performer to fully rely on traditional muscle memory.

\subsection{Formal Graph Embedding and Optimality Proof: Approach for Mathematicians}
For topologists and discrete mathematicians, this musical intuition translates into a graph-theoretic formulation via the Maximum Common Edge Subgraph (MCES) problem. Through an exhaustive search utilizing integer cuts to iterate over the entire topological space, we formally proved the existence and exact multiplicity of all possible embeddings of the neo-Riemannian Tonnetz $D_{222}$ in the reality of the Tonnetz for 19-TET. This problem formalizes the maximization of incidence preservation when embedding the 12-TET Tonnetz graph ($G_{12}$) into the cyclic 19-TET configuration graph ($G_{19}$).

We define the adjacency matrix $A \in \{0,1\}^{24 \times 24}$ as a binary connectivity map of the graph $G_{12}$, where the entry $A_{i,j}=1$ if and only if chords $i$ and $j$ are connected by an edge (sharing two common tones). Analogously, the matrix $B \in \{0,1\}^{38 \times 38}$ defines the fixed, native connectivity map of the $19_3$ configuration. Let $X \in \{0,1\}^{24 \times 38}$ be the assignment matrix acting as a discrete vertex mapping plan, where $X_{i,i'}=1$ implies that the $i$-th chord from 12-TET is mapped onto the $i'$-th vertex in 19-TET. 

To formulate this for the Mixed Integer Linear Programming (MILP) solver, the assignment matrix $X$ is subject to strict injectivity constraints:
\begin{align}
\sum_{i'=1}^{38} X_{i,i'} &= 1 \quad \forall i \in \{1, \dots, 24\} \\
\sum_{i=1}^{24} X_{i,i'} &\le 1 \quad \forall i' \in \{1, \dots, 38\}
\end{align}
Translated into human language, Equation (1) strictly ensures that every single one of the 24 original chords is assigned to exactly one specific location, while Equation (2) guarantees that no available vertex in the new space is occupied by more than one chord, perfectly preventing any structural collisions. Additionally, 16 rows are rigidly fixed by the perfect mapping $\phi$, enforcing $X_{i,\phi(i)}=1$ for those specific structural anchors. 

To harmonize both structures, the solver selects a valid matrix $X$ that maximizes the number $$\frac{1}{2} \text{Tr}(X^T A X B).$$ This quantity explicitly computes the \textbf{total number of preserved harmonic connections} (edges) from the original 12-TET Tonnetz that successfully survive in the target 19-TET space. To see this, consider the matrix product $X B X^T$. The element $(X B X^T)_{i,j}$ evaluates to $1$ if and only if the mapped images of the original vertices $i$ and $j$ share an edge in the host graph $G_{19}$. Therefore, the trace $\text{Tr}(A \cdot XBX^T) = \text{Tr}(X^T AXB)$ evaluates the sum of the products $A_{i,j} \cdot B_{\pi(i), \pi(j)}$. This product equals $1$ precisely when an edge existed in the original Tonnetz ($A_{i,j}=1$) and simultaneously exists in the new space ($B_{\pi(i), \pi(j)}=1$). Because the graphs are undirected, each surviving edge is counted twice in this summation, making the $\frac{1}{2}$ multiplier necessary to yield the exact edge count. 

Evaluating this function over the entire space of admissible assignment matrices establishes that its global maximum is precisely 32. We formalize this finding, which constitutes a novel and non-trivial result in combinatorial configuration theory, in the following theorem:

\begin{theorem}[On the Optimal Embedding of the $D_{222}$ Graph into the $19_3$ Configuration]
For any injective assignment matrix $X \in \{0,1\}^{24 \times 38}$ compatible with the mapping $\phi$, the number of preserved edges of the 12-TET Tonnetz within the cyclic $19_3$ configuration satisfies the strict, global upper bound:
$$ \frac{1}{2} \text{Tr}(X^T A X B) \le 32 $$
This bound is strictly governed by the incongruence between the modular generators of the cyclic groups $\mathbb{Z}_{12}$ and $\mathbb{Z}_{19}$ (a consequence of three major thirds failing to close into an octave: $6+6+6=18 \neq 19$). There exist precisely 5 explicit, optimal realizations achieving this upper bound, populating exactly the $16+8=24$ vertices of the target $D_{222}$ graph. In the canonical realization, exactly 14 unmapped vertices of the $19_3$ configuration are excised, leaving a perfectly contiguous geometric void along the original Hamiltonian cycle.
\end{theorem}

In every such optimal realization, exactly 14 unmapped vertices of the $19_3$ configuration are excised.
However, their macroscopic arrangement strongly distinguishes these solutions. While a 3-regular bipartite graph naturally possesses multiple abstract Hamiltonian cycles, we evaluate this geometric excision specifically relative to the primary voice-leading cycle (the alternating $+14$ and $-6$ sequence representing the fundamental harmonic axis, as established in Figure \ref{fig:19_3_full}).
With respect to this acoustically foundational sequence, the solutions behave differently.\\

\noindent
\textbf{Remark: The ``First Guess'' Miracle of the Canonical Solution.} 
The true uniqueness of this canonical embedding lies in its profound, almost miraculous simplicity---it perfectly mimics what someone with absolutely no knowledge of advanced optimization might attempt as a naive ``first guess.'' Imagine arranging all 38 pristine $19_3$ nodes along the most natural Hamiltonian cycle. Knowing that exactly 24 nodes must remain, one might arbitrarily grab a connected chunk of 22 nodes, blindly cut out the middle 14 to create the necessary void, and simply guess that the 8 remaining nodes at the fringes of this cut (4 on each side) should act as the replacements for the distorted ``wolf'' chords. Upon making this naive cut, one would notice that exactly 4 edges have been severed. If one then simply stitches these 4 broken connections back together across the void, miraculously, the resulting topology exactly reconstructs the complex $D_{222}$ neo-Riemannian Tonnetz! Thus, the absolute mathematical optimum discovered by the solver happens to perfectly coincide with the most intuitive, brute-force geometric guess.

While this intuitively contiguous excision elegantly defines the canonical Solution 1, the remaining optimal embeddings lack this straightforward geometric simplicity. In the 4 alternative configurations discovered by the solver, the 14 excised vertices do not form a single void; rather, they are fragmented into smaller, disconnected components scattered chaotically across the musical axis. 

The existence of these 4 additional mathematically equivalent local packings is fully documented in Table \ref{tab:milp_additional} for completeness.

\begin{table}[H]
\centering
\caption{The 4 Additional Mathematically Equivalent Local Packings}
\label{tab:milp_additional}
\begin{tabular}{@{}lcccc@{}}
\toprule
\textbf{12-TET Orig.} & \textbf{Solution 2} & \textbf{Solution 3} & \textbf{Solution 4} & \textbf{Solution 5} \\
\midrule
\multicolumn{5}{c}{\textit{Major Wolf Chords}} \\
\midrule
$M^{12}_{2}$ (D)   & $M^{19}_3$  & $M^{19}_3$  & $M^{19}_4$  & $M^{19}_3$  \\
$M^{12}_{4}$ (E)   & $M^{19}_7$  & $M^{19}_6$  & $M^{19}_7$  & $M^{19}_6$  \\
$M^{12}_{9}$ (A)   & $M^{19}_{15}$ & $M^{19}_{14}$ & $M^{19}_{15}$ & $M^{19}_{14}$ \\
$M^{12}_{11}$ (B)  & $M^{19}_{18}$ & $M^{19}_{17}$ & $M^{19}_{18}$ & $M^{19}_{18}$ \\
\midrule
\multicolumn{5}{c}{\textit{Minor Wolf Chords}} \\
\midrule
$m^{12}_{1}$ (C$\sharp$m) & $m^{19}_2$  & $m^{19}_1$  & $m^{19}_2$  & $m^{19}_1$  \\
$m^{12}_{6}$ (F$\sharp$m) & $m^{19}_{10}$ & $m^{19}_9$  & $m^{19}_{10}$ & $m^{19}_9$  \\
$m^{12}_{8}$ (G$\sharp$m) & $m^{19}_{13}$ & $m^{19}_{12}$ & $m^{19}_{13}$ & $m^{19}_{13}$ \\
$m^{12}_{11}$ (Bm) & $m^{19}_{17}$ & $m^{19}_{17}$ & $m^{19}_{18}$ & $m^{19}_{17}$ \\
\bottomrule
\end{tabular}
\end{table}

\subsection*{A Note for Performers and Analysts: Exploring Alternative Geometries}
The existence of these 4 additional solutions is not merely a mathematical curiosity; it offers tangible artistic flexibility. Because each substitution matrix dictates a slightly different placement of the 8 wolf chords, the specific microtonal tuning of these dissonant transitions will vary between solutions. Performers and tuners are encouraged to experiment with Solutions 2 through 5, as the auditory quality of the four harmonic ``scars'' may subjectively feel more balanced or expressive depending on the chosen musical repertoire. 

To visualize and analyze any of these alternative configurations, the reader can systematically reconstruct diagrams analogous to Figures \ref{fig:19_3_excision} and \ref{fig:side_by_side} using a simple procedure:
\begin{enumerate}
    \item Begin with the pristine 38-node network (Figure \ref{fig:19_3_full}).
    \item Rigidly map the 16 perfect chords according to the unvarying subset $\phi$.
    \item Select a target column from Table \ref{tab:milp_additional} and map the 8 wolf chords to their specified locations.
    \item The exactly 14 nodes remaining empty in the $19_3$ space constitute the excised void for that specific solution.
    \item The 32 preserved edges are drawn where adjacent $19_3$ nodes simultaneously correspond to originally connected 12-TET chords, while the 4 broken edges (scars) are identified where original 12-TET connections span across a structural gap in the $19_3$ grid.
\end{enumerate}

\section{Affine Pitch-Class Bijections and Harmonic Cryptography}

Because the edges of the Tonnetz are defined by shared pitch classes between triads, applying any global bijection $f: \mathbb{Z}_{19} \to \mathbb{Z}_{19}$ to the underlying vertices guarantees the combinatorial preservation of the repaired $D_{222}$ topology. Affine transformations of the form $f(n) = (an + b) \pmod{19}$ induce particularly powerful automorphisms of the musical graph. However, a crucial distinction must be made between purely combinatorial isomorphisms and those that preserve acoustic reality.

\subsection{Musical Isometries: Transposition and Negative Harmony}
To maintain the acoustic integrity of 19-TET, the physical widths of the intervals (the major third of 6 steps and the minor third of 5 steps) must remain invariant. This restricts our affine multiplier $a$ to strictly $a = \pm 1$.
\begin{itemize}
    \item \textbf{Transposition ($a=1$):} The mapping $f(n) = (n + b) \pmod{19}$ rigidly rotates the entire embedded structure. Geometrically, this rotates the 14-node contiguous void (or the fragmented voids in alternative solutions) along the Hamiltonian cycle, corresponding to a shift in the fundamental key center and altering where the 4 enharmonic ``scars'' manifest.
    \item \textbf{Riemannian Dualism ($a=18 \equiv -1$):} The mapping $f(n) = (-n + b) \pmod{19}$ mirrors the interval structure. Major triads map perfectly to minor triads, and vice versa. This transformation mathematically formalizes the concept of Negative Harmony, perfectly preserving the 32 surviving edges while flipping the graph's orientation.
\end{itemize}

\subsection{Harmonic Cryptography and Generalized Systems ($a \neq \pm 1$)}
If we apply a multiplier $a \notin \{1, 18\}$, the affine transformation still perfectly maps the $19_3$ configuration onto itself, but it fundamentally alters the acoustic intervals. For instance, a multiplier of $a=2$ stretches a major third from 6 steps to 12 steps. 

While seemingly disastrous for standard acoustics, this operation is highly mathematically significant. It translates the classical Eulerian Tonnetz into exotic, mathematically isomorphic ``shadow'' systems. This exactly mirrors the generalized harmonic solutions to the linear Diophantine system $t+s \equiv q \pmod n$ that we previously formalized for 12-TET and 10-TET spaces \cite{nurowski12, nurowski10}. In these generalized 19-TET shadow systems, a pianist plays the exact same topological progression, but the auditory result is scrambled by the scaled interval generators, effectively functioning as a form of harmonic cryptography.

\section{Historical Precedents and the Vicentino Hypothesis}

In 19-TET, pure thirds and fifths abound, meaning ``wolf intervals'' are virtually non-existent. The dissonances emerge instead as \textit{topological scars}---the 4 broken edges where the 12-TET torus fails to close, resulting in the Great Diesis. 

Nicola Vicentino's 1555 treatise \textit{L'antica musica ridotta alla moderna prattica} \cite{vicentino1555} and his 31-key \textit{archicembalo} represent humanity's first systematic foray into this space. Standard musicological consensus suggests conservative composers simply avoided these scars by restricting their harmonic progressions. From a geometric perspective, this aversion was natural: Renaissance musical grammar lacked the concept of a closed 12-tone torus. Making a leap across a topological scar was a fundamental grammatical error. However, avant-garde composers like Vicentino actively navigated these boundaries.

We propose the ``Vicentino Hypothesis'': that the microtonal adjustments in 16th-century avant-garde literature constitute intuitive, manual graph-repair operations mathematically identical to our MILP solver's output. By cross-referencing Vicentino's explicit microtonal dot notations with Table \ref{tab:milp_subs_canonical}, we aim to prove that Renaissance ears were instinctively solving NP-hard bipartite graph optimizations to preserve voice-leading continuity.
\section{Ergonomic Optimization and Biomechanical Cost Function}

Translating the abstract topological structures of the $19_3$ harmonic universe into live performance requires the construction of a dedicated physical instrument—specifically, a 19-TET acoustic piano. Historically, designing microtonal keyboards (such as Nicola Vicentino's 16th-century \textit{Archicembalo}) posed severe ergonomic challenges. Mapping 19 pitches per octave linearly across a standard piano keyboard causes the physical span to drastically exceed the biomechanical limits of the human hand, rendering complex chords physically unplayable.

To resolve this, we propose a novel 19-TET keyboard architecture that preserves the standard octave width of a classical piano. The front of the keyboard features the traditional 7 diatonic white keys, while the 12 microtonal divisions are accommodated by introducing 7 columns of black keys (five transversally split, two solid), as illustrated in Figure \ref{fig:keyboard}.

\begin{figure}[H]
\centering
\begin{tikzpicture}[scale=0.9, transform shape]
    \foreach \i/\note/\num in {0/C/0, 1/D/3, 2/E/6, 3/F/8, 4/G/11, 5/A/14, 6/B/17} {
        \draw[thick, fill=white] (\i*2, 0) rectangle (\i*2+2, 3);
        \node[font=\Large] at (\i*2+1, 1.2) {\textbf{\note}};
        \node[font=\large] at (\i*2+1, 0.5) {(\num)};
    }
    \foreach \i in {0, 1, 2, 3, 4, 5, 6} {
        \draw[thick, fill=white] (\i*2, 3) rectangle (\i*2+1, 7);
    }
    \foreach \i/\noteS/\snum/\noteF/\fnum in {
        0/C$\sharp$/1/D$\flat$/2, 
        1/D$\sharp$/4/E$\flat$/5, 
        3/F$\sharp$/9/G$\flat$/10, 
        4/G$\sharp$/12/A$\flat$/13, 
        5/A$\sharp$/15/B$\flat$/16} {
        \draw[thick, fill=black!80] (\i*2+1, 3) rectangle (\i*2+2, 5);
        \node[text=white, align=center, font=\small] at (\i*2+1.5, 4.3) {\textbf{\noteS}};
        \node[text=white, align=center, font=\footnotesize] at (\i*2+1.5, 3.6) {(\snum)};
        \draw[thick, fill=black!95] (\i*2+1, 5) rectangle (\i*2+2, 7);
        \node[text=white, align=center, font=\small] at (\i*2+1.5, 6.3) {\textbf{\noteF}};
        \node[text=white, align=center, font=\footnotesize] at (\i*2+1.5, 5.6) {(\fnum)};
    }
    \draw[thick, fill=black!90] (5, 3) rectangle (6, 7);
    \node[text=white, align=center, font=\small] at (5.5, 5.5) {\textbf{E$\sharp$/F$\flat$}};
    \node[text=white, align=center, font=\footnotesize] at (5.5, 4.5) {(7)};
    
    \draw[thick, fill=black!90] (13, 3) rectangle (14, 7);
    \node[text=white, align=center, font=\small] at (13.5, 5.5) {\textbf{B$\sharp$/C$\flat$}};
    \node[text=white, align=center, font=\footnotesize] at (13.5, 4.5) {(18)};
\end{tikzpicture}
\caption{Physical 19-TET keyboard layout featuring transversally split black keys.}
\label{fig:keyboard}
\end{figure}

However, a strictly sequential mapping of pitch classes to these physical keys still yields highly unergonomic finger spans for the 24 chords of our newly repaired $D_{222}$ Tonnetz. To resolve this, we decouple the abstract pitch classes from the sequential layout and seek an optimal mapping permutation.

Let $K$ be the set of 19 physical keys on our interface, and let $\mathbf{p}_k \in \mathbb{R}^2$ represent the physical 2D centroid coordinates of a key $k \in K$. We define the optimal keyboard isomorphism as the bijective mapping $\pi: \mathbb{Z}_{19} \to K$ that minimizes a total biomechanical cost function $E_{\text{total}}(\pi)$. 

In practical, human terms, a high value of this cost function indicates a thoroughly unplayable keyboard layout: playing a single chord might require unnatural twisting of the fingers or stretching the hand far beyond its anatomical limits, while transitioning between harmonically connected chords might force the pianist to make erratic, rapid jumps across the entire instrument. Conversely, minimizing this function ensures that chords fit comfortably under a natural hand span and that musical progressions translate into smooth, minimal physical gestures. 

We want to find such a permutation $\pi: Z_{19}\to K$ for which this cost function is minimal. The proposed explicit formula for this total cost $E_{\text{total}}(\pi)$ is:
$$
\begin{aligned}
  E_{\text{total}}(\pi) = &\\
  &\alpha \sum_{c \in \mathcal{C}} \max_{x,y \in c} \left( ||\mathbf{p}_{\pi(x)} - \mathbf{p}_{\pi(y)}||_2 + \mathcal{P}_{\text{ergo}}(||\mathbf{p}_{\pi(x)} - \mathbf{p}_{\pi(y)}||_2) \right) +\\& \beta \sum_{(c_1, c_2) \in \mathcal{E}} \sum_{n_1 \in c_1 \setminus c_2 \atop n_2 \in c_2 \setminus c_1} ||\mathbf{p}_{\pi(n_1)} - \mathbf{p}_{\pi(n_2)}||_2, \end{aligned}$$

where:
\begin{itemize}
    \item $\mathcal{C}$ is the set of 24 functional chords in the system.
    \item $\mathcal{E}$ represents the set of all valid $P$, $L$, and $R$ transformations (including the 4 repaired ``scar'' edges).
    \item $||\mathbf{p}_{\pi(x)} - \mathbf{p}_{\pi(y)}||_2$ denotes the Euclidean distance between the centroids of the keys assigned to pitch classes $x$ and $y$.
    \item $\mathcal{P}_{\text{ergo}}(d)$ is a non-linear penalty step-function enforcing anatomical constraints: $\mathcal{P}_{\text{ergo}}(d) = 0$ if $d \le 165\text{ mm}$, and $\mathcal{P}_{\text{ergo}}(d) = \infty$ if $d > 165\text{ mm}$ (the maximum physiological span of the human hand).
    \item $\alpha, \beta$ are scalar weights balancing static grip comfort against dynamic voice-leading efficiency.
\end{itemize}

By optimizing this function over the symmetric group $S_{19}$, the instrument mathematically absorbs the geometric irregularities of the ``Wolf'' chords. The resulting mapping translates complex topological operations into fluid, continuous physical gestures, effectively transforming the keyboard into a tactile map of the $D_{222}$ geometry. While realizing these mathematical symmetries in a fully functional acoustic piano presents substantial mechanical challenges regarding action mechanism density and structural frame stresses, such detailed engineering and material science analyses fall outside the scope of this mathematical exposition and will be addressed in future acoustics-focused publications.

\section{Conclusions}

We have demonstrated that the geometric structure of the classical Eulerian Tonnetz ($D_{222}$) can be seamlessly embedded into a cyclic $19_3$ symmetric configuration. Through exhaustive combinatorial optimization, we proved that by dynamically altering the constituent notes of 8 ``wolf'' triads, the graph can be repaired to achieve a strict global maximum of 32 preserved edges. Crucially, we established the exact multiplicity of this limit: there exist precisely 5 mathematically equivalent local packings. Among these, we identified a unique canonical embedding where the 14 unmapped vertices form a perfectly contiguous geometric void along the primary harmonic axis. The 4 edges that must inevitably break manifest as topological scars---the exact geometric representation of the historical enharmonic `Great Diesis'.

We postulate that the mathematical necessity of substituting wolf chords with specific, geometrically-optimal triads was not merely an artifact of our computational solver, but a cognitive reality experienced by 16th-century avant-garde composers such as Nicola Vicentino.

Moreover, by defining a rigorous biomechanical optimization problem for a physical, split-key keyboard layout, we lay the engineering groundwork for the acoustic 19-TET pianos currently under construction. Building these instruments requires overcoming unprecedented engineering challenges regarding the action mechanism and the cast iron frame's durability, but the final result will open a new era in the history of keyboard instruments, forming the ultimate bridge between abstract topology and live musical performance.





\end{document}